\documentclass[graybox]{svmult}

\usepackage{mathptmx}       % selects Times Roman as basic font
\usepackage{helvet}         % selects Helvetica as sans-serif font
\usepackage{courier}        % selects Courier as typewriter font
\usepackage{type1cm}        % activate if the above 3 fonts are
\usepackage{makeidx}         % allows index generation
\usepackage{graphicx}        % standard LaTeX graphics tool
\usepackage{multicol}        % used for the two-column index
\usepackage[bottom]{footmisc}% places footnotes at page bottom

\usepackage{amssymb}

\usepackage{amsmath}

\newtheorem{examples}[example]{Examples}
\newtheorem{assumption}[note]{Assumption}

\DeclareMathOperator{\supp}{supp}

\makeindex             % used for the subject index
\begin{document}

\title*{Elementary pathwise methods  for nonlinear parabolic and transport type SPDE with fractal noise}
\titlerunning{Elementary pathwise methods  for SPDEs with fractal noise} %for an abbreviated version of your contribution title if the original one is too long
\author{Michael Hinz \and Elena Issoglio \and Martina Z\"ahle}
\authorrunning{Michael Hinz et al.} %for an abbreviated version of
% your contribution title if the original one is too long
\institute{Michael Hinz  \at Department of Mathematics
Bielefeld University,
Postfach 100131, 33501 Bielefeld, Germany \email{MHinz@math.uni-bielefeld.de}
\and Elena Issoglio \at Department of Mathematics, King's College London, Strand, London, WC2R 2LS \email{Elena.Issoglio@kcl.ac.uk}
\and Martina Z\"ahle \at Institute of Mathematics
Friedrich Schiller University Jena
D-07737 Jena, Germany \email{Martina.Zaehle@uni-jena.de}}
%
% Use the package "url.sty" to avoid
% problems with special characters
% used in your e-mail or web address
%

\maketitle
\abstract{We survey some of our recent results on existence, uniqueness and regularity of function solutions to parabolic and transport type partial differential equations driven by non-differentiable noises. When applied pathwise to random situations, they provide corresponding statements for stochastic partial differential equations driven by fractional noises of sufficiently high regularity order. The approach is based on semigroup theory.}

%\tableofcontents

\section{Introduction}\label{S:Intro}

In this survey we list several of our recent results on existence, uniqueness and regularity of function solutions to linear and nonlinear parabolic stochastic partial differential equations such as abstract stochastic heat equations \cite{HZ09a, HZ09b, HZ12}, stochastic transport-diffusion equations \cite{I} and stochastic Burger's system, \cite{H09}. Our approach combines semigroup theory, \cite{P83, T78, Y80}, and fractional calculus, \cite{He, P, SKM93}. This leads to an elementary and easily accessible formulation in the sense that more sophisticated techniques such as rough path theory \cite{FV, L94, L98} are avoided and we obtain explicit formulas in terms of the semigroup. The basic idea of the studies surveyed here was to formulate a framework for stochastic partial differential equations using analogs of the pathwise techniques \cite{KZ, Z98, Z01, Z} previously employed by the third named author to solve stochastic differential equations.

General information on stochastic partial differential equations can be found in \cite{DPZ92, HOUZ96, PR07, W86}, results close to our discussion of parabolic equations are for instance \cite{GLT06, HN08, MN03, TTV03}. A transport type equation was investigated in \cite{RT07}, and some results on stochastic Burgers equation can be found in \cite{BCJL94, DDT94, DG95, GN99}. Of course there are many more highly valuable references on these topics.

Classical stochastic calculus allows to integrate predictable processes against semimartingale integrators. In particular it can be used to study stochastic differential equations with respect to a given semimartingale process. From a practical viewpoint both semimartingale properties of the integrator and predictability of the integrand may be too restrictive in some situations. If the integrator is Gaussian we may use methods from Malliavin calculus to define stochastic integrals for nonanticipating integrands, \cite{N}. Alternatively, if almost surely both the integrand and the integrator are of sufficiently high regularity (for instance in the sense of H\"older continuity, $p$-Variation or fractional differentiability) then this regularity can be used to define stochastic integrals in a pathwise sense. In this case they are of Stieltjes type. By now the most popular approach to this idea probably is Young integration, \cite{Y36}, which later inspired the development of rough path theory, \cite{L94, L98}. 
Another way to go, although not entirely pathwise, is to use stochastic calculus via regularization, \cite{RV93}. Yet another technique was introduced in \cite{Z98, Z01} and is based on fractional calculus.

Let $I=(a,b)$ be a bounded interval and $E$ be a Banach space. Given $\eta>0$ and a function $\varphi \in L_1((a,b),E) $, consider the (forward and backward) \emph{Riemann-Liouville fractional integrals of order $\eta$} by
\[\mathrm{I}_{a+}^\eta \varphi(t):=\frac{1}{\Gamma(\eta)}\int_a^t\frac{\varphi(\tau)}{(t-\tau)^{1-\eta}}\D\tau\]
and
\[\mathrm{I}_{b-}^\eta\varphi(t):=\frac{(-1)^{-\eta}}{\Gamma(\eta)}\int_t^b \frac{\varphi(\tau)}{(\tau-t)^{1-\eta}}\D\tau\; . \]
Here for $\eta>0$ the powers are understood as usual in the sense of choosing the main branch of the analytic function $\zeta^\eta$, $\zeta\in \mathbb{C}$, with the cut along the positive half axis, in particular, $(-1)^\eta=\E^{i\eta\pi}$. Here and in the following the integrals are understood in the Bochner sense. Let $\mathrm{I}^\eta_{a+}(L_p((a,b),E)) $ denote the space of functions $f=\mathrm{I}_{a+}^\eta\varphi $ with $\varphi \in L_p((a,b),E) $, similarly $\mathrm{I}_{b-}^\eta(L_p((a,b),E))$. For $0<\eta< 1 $ and functions $f \in \mathrm{I}^\eta_{a+}(L_p((a,b),E)) $, respectively $f \in \mathrm{I}^\eta_{b-}(L_p((a,b),E)) $, consider the  \emph{left-sided Weyl-Marchaud fractional derivatives of order $\eta $},
\[\mathrm{D}_{a+}^\eta f(t):=\frac{\mathbf{1}_{(a,b)}(t)}{\Gamma(1-\eta)}\left(\frac{f(t)}{(t-a)^\eta}+ \eta\int_a^t\frac{f(t)-f(\tau)}{(t-\tau)^{\eta+1}}\D\tau\right)  \]
and the \emph{right-sided Weyl-Marchaud fractional derivatives of order $\eta $},
\[\mathrm{D}_{b-}^\eta f(t):=\frac{(-1)^\eta\mathbf{1}_{(a,b)}(t)}{\Gamma(1-\eta)}\left(\frac{f(t)}{(b-t)^\eta}+ \eta\int_t^b\frac{f(t)-f(\tau)}{(\tau-t)^{\eta+1}}\D\tau\right)\;,\]
the convergence of the principal values being pointwise almost everywhere if $p=1$ and in $L_p((a,b),E)$ if $p > 1$. Under these assumptions $\mathrm{I}_{a+}^\eta \mathrm{D}_{a+}^\eta f=f $ in $L_p((a,b),E)$, while $\mathrm{D}_{a+}^\eta \mathrm{I}_{a+}^\eta \varphi=\varphi $ is true for any $\varphi\in L_1((a,b),E) $.
In the case $\eta=1$ set $\mathrm{D}_{a+}^1f=\D f/\D t $ and $\mathrm{D}_{b-}^1 f=-\D f/\D t $ and in the case $\eta=0$ define $\mathrm{D}_{a+}^0$ and $\mathrm{D}_{b-}^0$ to be the identity.
See for instance \cite{HZ09a, SKM93}.

For a moment assume that $E=\mathbb{R}$ and consider real valued functions $f$ and $g$ on $(a,b)$ such that the limits $f(a+)$, $g(a+)$ and $g(b-)$ exist. Consider the regulated functions
\begin{equation}\label{E:corrections}
f_{a+}(t):=\mathbf{1}_{(a,b)}(t)(f(t)-f(a+))\ \ \text{ and } \ \ g_{b-}(t):=\mathbf{1}_{(a,b)}(t)(g(t)-g(b-))\; .
\end{equation}
In \cite{Z98} it had been shown that if $f_{a+}\in \mathrm{I}_{a+}^\eta(L_p(a,b))$ and $g_{b-}\in \mathrm{I}_{b-}^{1-\eta}(L_q(a,b))$ for some $1/p+1/q\leq 1$ and $0\leq \eta\leq 1$ then the integral
\begin{equation}\label{E:forwardint}
\int_a^b f(s)\D g(s):=(-1)^\eta \int_a^b \mathrm{D}_{a+}^\eta f_{a+}(s) \mathrm{D}_{b-}^{1-\eta}g_{b-}(s)\D s +f(a+)(g(b-)-g(a+))
\end{equation}
is well defined, that is, the value of the right-hand side in (\ref{E:forwardint}) is a real number that is independent of the particular choice of $\eta$. Moreover, if $f$ and $g$ are sufficiently regular such that both (\ref{E:forwardint}) and the Lebesgue-Stieltjes integral $(LS)\int_a^bf\D g$ exist, then they agree. For instance, if $f$ and $g$ are H\"older continuous and the sum of their H\"older orders is greater than one, then (\ref{E:forwardint}) exists and equals the Riemann-Stieltjes integral $(RS)\int_a^bf\D g$. If $f$ and $g$ satisfy the above conditions with $0\leq \eta< 1/p$ then the correction terms in (\ref{E:forwardint}) may be dropped, more precisely, we have
\[
\int_a^bf(s)\D g(s)=(-1)^\eta \int_a^b \mathrm{D}_{a+}^\eta f(s) \mathrm{D}_{b-}^{1-\eta}g_{b-}(s)\D s\; .
\]
See \cite{Z98, Z01} for details. Integrals of type (\ref{E:forwardint}) may for instance be used to investigate differential equations of form
\begin{equation}\label{E:ode}
\begin{cases}
\D x(t) = a(x(t),t)\D z(t)+b(x(t),t)\D t  \\
x(0)=x_0\;,
\end{cases}
\end{equation}
where $z$ is a non-differentiable function that is H\"older continuous of order greater than $1/2$ and $a$ and $b$ are coefficients that satisfy certain growth and smoothness assumptions. Equation (\ref{E:ode}) is made precise by saying that $x=(x(t))_{t\geq 0}$ solves (\ref{E:ode}) if
\[x(t)=x_0+\int_0^ta(x(s),s)\D z(s)+\int_0^tb(x(s),s)\D s\]
for any $t>0$, where the first integral is defined as in  (\ref{E:forwardint}). As usual, the existence and uniqueness of solutions $x$ to (\ref{E:ode}) is obtained by combining a priori estimates on the integral operator $x\mapsto \int_0^\cdot a(x(s),s)\D z(s)$ in suitable function spaces and fixed point arguments, \cite{KZ, Z01}. If typical realizations of suitable random processes are used in place of $z$, such as for instance the paths of a fractional Brownian motion $B^H$ with Hurst parameter $H>1/2$, this yields a stochastic differential equation in the pathwise sense.

Equation (\ref{E:ode}) is an evolution problem subject to a perturbation $z$. Also parabolic partial differential equations of form
\[\begin{cases}
\frac{\partial u}{\partial t}(t) =-\mathrm{A}u(t)\\
u(0)=u_0\end{cases}\]
are commonly viewed as evolution problems, now of course in abstract (Banach or Hilbert) spaces, \cite{P83}, and their behaviour is completely governed by a related semigroup $(\mathrm{T}(t))_{t\geq 0}$ of evolution operators, that is, the solution $u$ to the Cauchy problem will be of the form $u(t)=\mathrm{T}(t)u_0$, $t>0$. We will use an analog of (\ref{E:forwardint}) to incorporate a noise signal $z$ into the equation. A simple linear multiplicative perturbation would for instance lead to a Cauchy problem of form
\[\begin{cases}
\frac{\partial u}{\partial t}(t) =-\mathrm{A}u(t)+u(t)\cdot \dot{z}(t)\\
u(0)=u_0\;.\end{cases}\]
If the noise $z$ is random, this yields again a pathwise technique, now for stochastic partial differential equations. It allows to investigate problems perturbed by signals that lack semimartingale properties but have sufficiently high regularity in terms of H\"older and Sobolev norms.

To consider a version of (\ref{E:forwardint}) for vector valued functions let $E$ and $F$ be se\-pa\-ra\-ble Banach spaces and let $L(E,F)$ denote the space of bounded linear operators from $E$ into $F$. Given $0\leq \eta\leq 1$, an $E$-valued function $z$ on $(a,b)$ and an $L(E,F)$-valued function $U$ on $(a,b)$ such that $\mathrm{D}_{b-}^{1-\eta}z_{b-}\in L_\infty((a,b),E)$ and $\mathrm{D}_{a+}^\eta U\in L_1((a,b),L(E,F))$, the integral
\begin{equation}\label{E:vectorforward}
\int_a^bU(s)\D z(s):=(-1)^\eta\int_a^b \mathrm{D}_{a+}^\eta U(s) \mathrm{D}_{b-}^{1-\eta}z_{b-}(s)\D s
\end{equation}
is well defined. More precisely, the right-hand side of (\ref{E:vectorforward}) is an element of $F$ and does not depend on the particular choice of $\eta$. The notation $z_{b-}$ is to be understood as in (\ref{E:corrections}).

\section{Semigroups and function spaces}\label{S:Prelim}

Let $(X, \mathcal{X}, \mu)$ be a $\sigma$-finite measure space and let $L_p(\mu)$, $1<p<\infty$ and $L_\infty(\mu)$ denote the spaces of (equivalence classes of) $p$-integrable and essentially bounded functions on $X$, respectively.

We assume that $\mathrm{T}=(\mathrm{T}(t))_{t\geq 0}$ is a symmetric strongly continuous semigroup on $L_2(\mu)$, that is, $\mathrm{T}(t+s)=\mathrm{T}(t)\circ \mathrm{T}(s)$, $\mathrm{T}(0)u=u$ and $\left\langle \mathrm{T}(t)u,v\right\rangle_{L_2(\mu)}=\left\langle u,\mathrm{T}(t)v\right\rangle_{L_2(\mu)}$ for any $s,t\geq 0$ and any $u,v\in L_2(\mu)$ and $\lim_{t\to 0}\left\|\mathrm{T}(t)u-u\right\|_{L_2(\mu)}=0$ for any $u\in L_2(\mu)$. We further assume that $(\mathrm{T}(t))_{t\geq 0}$ is Markovian, that is for any $t\geq 0$ and any $u\in L_2(\mu)$ with $0\leq u\leq 1$ $\mu$-a.e. we have $0\leq \mathrm{T}(t)u\leq 1$ $\mu$-a.e. In this case the semigroup $(\mathrm{T}(t))_{t\geq 0}$ is automatically contractive, $\left\|\mathrm{T}(t)u\right\|_{L_2(\mu)}\leq \left\|u\right\|_{L_2(\mu)}$, $t\geq 0$, $u\in L_2(\mu)$.

Let $-\mathrm{A}$ denote the infinitesimal $L_2(\mu)$-generator of $(\mathrm{T}(t))_{t\geq 0}$,
\[-\mathrm{A}u=\lim_{t\to 0} \frac{1}{t}(\mathrm{T}(t)u-u)\ \text{ strongly in $L_2(\mu)$}\]
for members $u$ of $\mathcal{D}(\mathrm{A})$, the dense subspace of $L_2(\mu)$ for whose members this limit exists. Both $\mathrm{A}$ and $\mathrm{T}(t)$ are non-negative definite self-adjoint operators on $L_2(\mu)$. In particular, the fractional powers $\mathrm{A}^\alpha$, $\alpha\geq 0$, of $\mathrm{A}$ can be defined in the usual way using the spectral representation, \cite{T78, Y80}.

For $0<\alpha<1$, we can characterize the domain $\mathcal{D}(\mathrm{A}^\alpha)$ of $\mathrm{A}^\alpha$ in terms of the semigroup: $u\in L_2(\mu)$ is in $\mathcal{D}(\mathrm{A}^ \alpha)$ if and only if
\begin{equation}\label{E:pospow}
\mathrm{A}^\alpha u=\lim_{\varepsilon\to 0}\frac{1}{\Gamma(-\alpha)}\int_\varepsilon^ \infty t^{-\alpha-1}(\mathrm{T}(t)-\mathrm{I})u \D t
\end{equation}
converges in $L_2(\mu)$, see e.g. \cite{BBW68}. This may be interpreted as a right-sided Weyl-Marchaud derivative $\mathrm{D}^\alpha_-$ of $t\mapsto \mathrm{T}(t)u$ at $t=0$, more precisely, $\mathrm{D}_-^ \alpha(\mathrm{T}(\cdot)u))(t)=(-1)^\alpha \mathrm{A}^\alpha \mathrm{T}(t)u$. See \cite{HZ09b} or \cite{SKM93}. Now let us temporarily assume that zero is not an eigenvalue of $\mathrm{A}$. Then the negative fractional powers $\mathrm{A}^{-\alpha}$, $\alpha>0$, can be expressed by
\begin{equation}\label{E:negpow}
\mathrm{A}^{-\alpha}u=\frac{1}{\Gamma(\alpha)}\int_0^\infty t^{\alpha-1}\mathrm{T}(t)u \D t\; ,
\end{equation}
what may be read as a right-sided Riemann-Liouville integral $\mathrm{I}_-^\alpha$ of order $\alpha>0$ of the function $t\mapsto \mathrm{T}(t)u$, i.e. $\mathrm{I}_-^\alpha(T(\cdot)u)(t)=(-1)^{-\alpha}\mathrm{A}^{-\alpha}\mathrm{T}(t)u$. Thus, for semigroups the language of traditional fractional calculus just leads to special cases of the usual functional calculus, cf.~\cite{SKM93,St70,Y80}.

The contractivity implies that $(\mathrm{T}(t))_{t\geq 0}$ is analytic on $L_2(\mu)$, cf. \cite{D89} or \cite{St70}, Chapter III. It also defines analytic semigroups on the spaces $L_p(\mu)$, $1\leq p<\infty$, see \cite{D89}, Theorem 1.4.1. or \cite{St70}, Chapter III. We use the same notation $\mathrm{T}=(\mathrm{T}(t))_{t\geq 0}$ for these semigroups but denote their $L_p(\mu)$-generators by $-\mathrm{A}_p$, such that $\mathrm{A}_2=\mathrm{A}$. In these cases (\ref{E:pospow}) and (\ref{E:negpow}) may be used to define their fractional powers, see  \cite{Y80}. Analyticity implies further useful properties: For any $u\in L_p(\mu)$, any $\alpha\geq 0$ and any $t>0$ we have
\begin{equation}\label{E:smooth}
\mathrm{T}(t)u\in\mathcal{D}(\mathrm{A}_p^\alpha)\; .
\end{equation}
The operators $\mathrm{T}(t)$ and $\mathrm{A}_p^\alpha$ commute on $\mathcal{D}(\mathrm{A}_p^\alpha)$. Given $\omega>0$, the bound
\begin{equation}\label{E:sgbound}
\left\|(\omega \mathrm{I}+\mathrm{A}_p)^\alpha \mathrm{T}(t)\right\|\leq c_\alpha \E^{\omega t} t^{-\alpha}
\end{equation}
holds for $t>0$ (in the operator norm on $L_p(\mu)$) and the continuity estimate
\begin{equation}\label{E:sgcont}
\left\|\mathrm{T}(t)u-u\right\|_{L_p(\mu)}\leq c_\alpha t^\alpha\left\|(\omega \mathrm{I}+\mathrm{A}_p)^\alpha u\right\|_{L_p(\mu)}+(1-\E^{-\omega t})\left\|u\right\|_{L_p(\mu)}
\end{equation}
is valid for $0\leq \alpha<1$, $u\in \mathcal{D}(\mathrm{A}_p^\alpha)$ and $t>0$. See \cite{P83}.

Given $\alpha_1,\alpha_2\geq 0$, we have $\mathrm{A}_p^{\alpha_1+\alpha_2}=\mathrm{A}_p^{\alpha_1}\mathrm{A}_p^{\alpha_2}$, $\mathrm{A}_p^{\alpha_1}\mathrm{A}_p^{-\alpha_1}=\mathrm{I}$ and $\mathrm{A}_p^{\alpha_1}:\mathcal{D}(\mathrm{A}_p^{\alpha_1+\alpha_2})\to \mathcal{D}(\mathrm{A}_p^{\alpha_2})$ is an isomorphism between these domains endowed with the graph norm. For $\sigma\geq 0$ we may regard the negative power
\[\mathrm{J}_p^\sigma(\mu):=(\mathrm{A}_p+ \mathrm{I})^{-\sigma/2}\; .\]
as a \emph{generalized Bessel potential operator} on $L_p(\mu)$. Set
\begin{equation}\label{E:revpotentials}
H_p^\sigma(\mu):=\mathrm{J}_p^\sigma(\mu)(L_p(\mu))\; ,
\end{equation}
$\sigma\geq 0$, equipped with the norms
\[\left\|u\right\|_{H_p^\sigma(\mu)}:=\left\|u\right\|_{L_p(\mu)}+ \left\|\mathrm{A}_p^{\sigma/2}u\right\|_{L_p(\mu)}\; .\]
Clearly $H^0_p(\mu)=L_p(\mu)$. If $p=2$ we write $H^\sigma(\mu)$ for $H_2^\sigma(\mu)$. Note that $\mathcal{D}((\mathrm{I}+\mathrm{A}_p)^\alpha)=\mathcal{D}(\mathrm{A}_p^\alpha)=H_p^{2\alpha}(\mu)$ and that potential operators $\mathrm{J}_p^\sigma(\mu)$, $\sigma\geq 0$ define isomorphic mappings from $H_p^{\alpha}(\mu)$ onto $H_p^{\alpha+\sigma}(\mu)$, $\alpha\geq 0$. Subspaces of essentially bounded functions will be denoted by
\[H_{\infty}^\sigma(\mu):=H^\sigma(\mu)\cap L_\infty(\mu)\; ,\]
normed by $\left\|\cdot\right\|_{H_{\infty}^\sigma(\mu)}:=\left\|\cdot\right\|_{H^\sigma(\mu)}+ \left\|\cdot\right\|_{L_\infty(\mu)}$. We write
\begin{equation}\label{E:revduals}
H_{p'}^{-\sigma}(\mu):=((H_p^\sigma(\mu))^\ast
\end{equation}
for the \emph{duals of the spaces} $H_p^\sigma(\mu)$, $1<p<\infty$, $\sigma\geq 0$, ${1}/{p}+{1}/{p'}=1$, equipped with the usual (operator) norm  $\left\|\cdot\right\|_{H_{p'}^{-\sigma}(\mu)}$.

If $X=\mathbb{R}^n$ and $\mu$ is the $n$-dimensional Lebesgue measure then the spaces $H_p^\sigma(\mu)$, $1<p<\infty$, $\sigma\geq 0$, coincide with potential spaces defined in terms of Fourier analysis,
\[H_p^\sigma(\mathbb{R}^n):=\left\lbrace \mathcal{S}'(\mathbb{R}^n): (1+|\xi|^2)^{\sigma/2}\hat{f})^\vee \in L_p(\mathbb{R}^n)\right\rbrace \; ,\]
$1<p<\infty$, $\sigma\in\mathbb{R}$. Here $f\mapsto \hat{f}$ and $f\mapsto \check{f}$ denote the Fourier transform and the inverse Fourier transform, and $\mathcal{S}'(\mathbb{R}^n)$ is the space of Schwartz distributions on $\mathbb{R}^n$. Given a smooth bounded domain $D\subset\mathbb{R}^n$ we also consider the spaces
\[\widetilde{H}^\sigma_p(D):=\left\lbrace f\in H_p^\sigma(\mathbb{R}^n): \supp f \subset \overline{D}\right\rbrace, \]
which are defined as subspaces of $H^\sigma_p(\mathbb{R}^n)$ for any $\sigma>-1/p$. We write $\widetilde{H}^\sigma(D)$ if $p=2$. The spaces $\widetilde{H}^\sigma_p(D)$ may be regarded as the potential spaces associated with the operator $\mathrm{A}_p$ given as the $L_p$-generator of the Dirichlet heat semigroup $(T^D(t))_{t\geq 0}$ for $D$. For $\alpha,\sigma\in\mathbb{R}$ with $-1/2<\alpha$ and $\alpha-\sigma<3/2$ the fractional power $\mathrm{A}^{\sigma/2}$ maps $\widetilde{H}^\alpha(D)$ isomorphically onto $\widetilde{H}^{\alpha-\sigma}$. If $0\leq \alpha<\frac32$ and $\alpha\neq\frac12$ then $\mathcal{D}(\mathrm{A}^{\alpha/2})=\widetilde{H}^\alpha(D)$. See \cite{T78}. The analyticity of $(\mathrm{T}^D(t))_{t\geq 0}$ also implies that for $-1/2<\alpha,\sigma,\sigma+\alpha<3/2$ the semigroup operators $\mathrm{T}^D(t)$ map $\widetilde{H}^\alpha(D)$ into $\widetilde{H}^{\alpha+\sigma}(D)$. In particular, given $f\in \widetilde{H}^\alpha(D)$ we will have $\supp \mathrm{T}^D(t)f\subset \overline{D}$.

In the following we will exclusively use subspaces consisting of real-valued functions and real valued dual elements respectively distributions. For simplicity we will not emphasize this fact by introducing new notation and therefore ask the reader to keep it in mind.

As we are going to investigate semilinear and transport equations, we need some preliminaries on composition and multiplication. Let $0\leq \sigma\leq 1$. If $F\in C(\mathbb{R})$ satisfies $F(0)=0$ and is Lipschitz, then we have
\[\left\|F(u)\right\|_{H^\sigma_\infty(\mu)}\leq c\left\|u\right\|_{H^\sigma_\infty(\mu)}\]
for any $u\in H^\sigma_\infty(\mu)$. If $F\in C^1(\mathbb{R})$ is such that $F(0)=0$ and its derivative $F'$ is bounded and Lipschitz, then
\[\left\|F(u)-F(v)\right\|_{H^\sigma_\infty(\mu)}\leq c\left\|u-v\right\|_{H^\sigma_\infty(\mu)} \left(\left\|v\right\|_{H^\sigma_\infty(\mu)}+1\right)\]
for any $u,v\in H^\sigma_\infty(\mu)$. Finally, if $F\in C^2(\mathbb{R})$ with $F(0)=0$ and bounded and Lipschitz second derivative $F''$, then
\[\left\| F(u_1)-F(v_1)-F(u_2)+F(v_2)\right\|_{H^\sigma_\infty(\mu)}\leq c\left(\left\|u_1-v_1-u_2+v_2\right\|_{H^\sigma_\infty(\mu)}+\left\|u_2- v_2\right\|_{H^\sigma_\infty(\mu)}\right)\]
for all $u_1,v_1,u_2,v_2\in H^\sigma_\infty(\mu)$ with $\|u_i\|_{H^\sigma_\infty(\mu)}\leq 1$ and $\|v_i\|_{H^\sigma_\infty(\mu)}\leq 1$ for $i=1,2$. These properties basically follow from the Markov property and the mean value theorem, see \cite[Proposition 3.1]{HZ12}. If $u,v\in H^\sigma_\infty(\mu)$, $0\leq \sigma\leq 1$, then again by the Markov property the pointwise product $uv$ is again in $H^\sigma_\infty(\mu)$ and
\[\left\|uv\right\|_{H^\sigma_\infty(\mu)}\leq c\left\|u\right\|_{H^\sigma_\infty(\mu)}\left\|v\right\|_{H^\sigma_\infty(\mu)}\; .\]
Given $u\in H^\sigma_\infty(\mu)$ and $z\in (H^\sigma_\infty)^\ast$, $0\leq \sigma\leq 1$, we can define the product $uz \in (H^\sigma_\infty(\mu))^\ast$ by
\[(uz)(v):=\left( z, uv\right), \ \ v\in H^\sigma_\infty(\mu)\; ,\]
where $\left(\cdot,\cdot\right)$ denotes the dual pairing. For $z\in H^{-\sigma}(\mu)$ we observe
\begin{equation}\label{E:basicproduct}
\left\|uz\right\|_{(H^\sigma_\infty(\mu))^\ast}\leq \left\|z\right\|_{H^{-\sigma}(\mu)}\left\|u\right\|_{H^\sigma_\infty(\mu)},
\end{equation}
note that $H^{-\sigma}(\mu)$ is a subspace of $(H^\sigma_\infty(\mu))^\ast$.

The semigroup $(\mathrm{T}(t))_{t\geq 0}$ is called \emph{(locally) ultracontractive with spectral dimension $d_S>0$} if there exist constants $c>0$ and $0<\omega\leq 1$ such that for any $t>0$ we have
\begin{equation}\label{E:ultra}
\left\|\mathrm{T}(t)\right\|_{L_2(\mu)\to L_\infty(\mu)}\leq ct^{-d_S/4}\E^{\omega t}.
\end{equation}
The estimate (\ref{E:ultra}) is equivalent to several functional inequalities of Nash and Sobolev type, see \cite{CKS87, C96, D89, V85}. If (\ref{E:ultra}) holds we can define $\mathrm{T}(t)z$ for $z\in (H^\sigma_\infty(\mu))^\ast$ by means of dual pairing,
\[(\mathrm{T}(t)z)(v):=\left( z, \mathrm{T}(t)v\right), \quad v\in L_2(\mu)\; ,\]
where we have implicitly used (\ref{E:smooth}). For $z\in H^{-\sigma}(\mu)$ we obtain
\begin{equation}\label{E:sgdualbound}
\left\|\mathrm{T}(t)z\right\|_{L_2(\mu)}\leq c\E^{\omega t} (t^{-\sigma/2}+t^{-d_S/4})\left\|z\right\|_{H^{-\sigma}(\mu)}
\end{equation}
by (\ref{E:sgbound}) and (\ref{E:ultra}).

\section{Integral operators}\label{S:Intops}

Using some of the facts from the preceding section allows to verify the existence of a version of (\ref{E:vectorforward}) that is suitable to solve related parabolic problems, \cite{HZ12}.

Let $t>0$ and assume that $u$ is a function on $(0,t)$ taking values in $H_\infty^\delta(\mu)$ for some $0<\delta<1$. If moreover $w\in H^{-\beta}(\mu)$ with $0<\beta\leq\delta$ and $G\in C(\mathbb{R})$ is Lipschitz with $G(0)=0$, then $G(u(\cdot))w$ is a function on $(0,t)$ taking its values in
$(H_\infty^\beta(\mu))^\ast$. By (\ref{E:sgdualbound})
\[s\mapsto U(t;s)w:=\mathrm{T}(t-s)G(u(s))w\ ,\ w\in H^{-\beta}(\mu),\]
is seen to define a function $s\mapsto U(t;s)$ that takes its values in $L(H^{-\beta}(\mu),H_\infty^\delta(\mu))$. If $s\mapsto u(s)$ is sufficiently regular,
\[\mathrm{D}_{0+}^\eta U(t;s):=\frac{\mathbf{1}_{(0,t)}(s)}{\Gamma(1-\eta)}\left( \frac{U(t;s)}{s^\eta}+\eta\int_0^s\frac{U(t;s)-U(t;\tau)}{(s-\tau)^{\eta+1}}\D\tau\right) \]
converges in an appropriate sense.

To make this more precise, we introduce some additional function spaces. Given a separable Banach space $E$ with norm $\left\|\cdot\right\|_E$ and a number $0<\eta<1$ let $W^\eta([0,t_0], E)$ denote the space of  $E$-valued functions $v$ on $[0,t_0]$ such that
\[\left\|v\right\|_{W^\eta([0,t_0],E)}:=\sup_{0\leq t\leq t_0}\left( \left\|v(t)\right\|_{E}+\int_0^t \frac{\left\|v(t)-v(\tau)\right\|_E}{(t-\tau)^{\eta+1}}\D\tau\right) <\infty\; .\]
Similarly, let $C^\eta([0,t_0],E)$, $0<\eta<1$, denote the space of $\eta$-H\"older continuous $E$-valued functions $v$ on $[0, t_0]$ such that
\[\left\|v\right\|_{C^\eta([0,t_0],E)}:=\sup_{0\leq t\leq t_0}\left\|v(t)\right\|_{E}+\sup_{0\leq\tau<t\leq t_0}\frac{\left\|v(t)-v(\tau)\right\|_E}{(t-\tau)^\gamma}<\infty\; .\]

\begin{lemma}\label{L:WeylMarchaud}
Let $0<\eta<1$, $t\in (0,t_0)$ and let $G\in C^2(\mathbb{R})$ with $G(0)=0$ and bounded and Lipschitz second derivative $G''$. If $0<\beta\leq \delta<1$, $u\in W^\eta([0,t], H_{\infty}^\delta(\mu))$ and
\[\delta\vee \frac{d_S}{2}<2-2\eta-\left(\beta\vee\frac{d_S}{2}\right)\; ,\]
then $\mathrm{D}_{0+}^\eta U(t;\cdot)$ converges in $L_1([0,t], L(H^{-\beta}(\mu), H_\infty^{\delta}(\mu)))$ and admits the following representation in terms of the semigroup:
\begin{multline}
\mathrm{D}_{0+}^\eta U(t;s)=\mathrm{D}_{0+}^\eta(\mathrm{T}(t-\cdot)G(u(\cdot))(s)\notag\\
=\mathbf{1}_{(0,t)}(s)\left\lbrace -\mathrm{A}^\eta \mathrm{T}(t-s)G(u(s))+ c_\eta \mathrm{T}(t-s)\int_s^\infty r^{-\eta-1}\mathrm{T}(r)G(u(s))\D r \right.\\
\left.+ c_\eta\int_0^s r^{-\eta-1}\mathrm{T}(r+t-s)[G(u(s))-G(u(s-r))]\D r\right\rbrace\; .
\end{multline}
Here $c_\eta=\eta \Gamma(1-\eta)^{-1}=-\Gamma(-\eta)^{-1}$.
\end{lemma}

Given $z \in C^{1-\alpha}([0,t_0], H^{-\beta}(\mu))$ and $\eta$ slightly bigger than $\alpha$ we may consider
\[\mathrm{D}_{t-}^{1-\eta}z_t(s):=\frac{(-1)^{1-\eta}\mathbf{1}_{(0,t)}(s)}{\Gamma(\eta)}\left( \frac{z(s)-z(t)}{(t-s)^{1-\eta}}+(1-\eta)\int_ s^t\frac{z(s)-z(\tau)}{(\tau-s)^{(1-\eta)+1}}\D\tau\right) \; ,\]
where $z_t(s):=\mathbf{1}_{(0,t)}(s)(z(s)-z(t))$. Then
\begin{equation}\label{E:w}
w(s):=\mathrm{D}_{t-}^{1-\eta}z_t(s)\ ,\ s\in [0,t]
\end{equation}
defines a function in $L_\infty([0,t], H^{-\beta}(\mu))$.

The next definition introduces an integral operator that may be seen as a version of (\ref{E:vectorforward}). Recall the notation $U(t;s)=\mathrm{T}(t-s)G(u(s))$.

\begin{definition}\label{D:integral}
Given $t\in [0,t_0]$, $0<\eta<1$ and sufficiently regular functions $u$ and $z$ on $[0,t]$, put
\begin{equation}\label{E:integral}
\int_0^t \mathrm{T}(t-s) G(u(s))\D z(s):=(-1)^\eta\int_0^t \mathrm{D}_{0+}^\eta U(t;s)\mathrm{D}^{1-\eta}_{t-} z_t(s)\D s\; .
\end{equation}
\end{definition}

This integral operator is well defined.

\begin{lemma}\label{L:existence}
Let $t$ and $\eta$ be as in Definition \ref{D:integral}. Assume $u$ is such that $\mathrm{D}_{0+}^\eta U(t;\cdot)\in L_1([0,t], L(H^{-\beta}(\mu), H_\infty^{\delta}(\mu)))$ and $z$ is such that $\mathrm{D}_{t-}^{1-\eta}z_t\in L_\infty([0,t], H^{-\beta}(\mu))$, where $0<\beta\leq\delta<1$.
Then the right-hand side of (\ref{E:integral}) exists as an element of $H_\infty^{\delta}(\mu)$ and is independent of the particular choice of $\eta$.
\end{lemma}

The following contraction property can be used to prove the existence and uniqueness of function solutions to Cauchy problems related to perturbed parabolic equations. To establish it we use equivalent norms on the space $W^\eta([0,t_0], E)$, $0<\eta<1$, given by
\[\left\|v\right\|^{(\varrho)}_{W^\eta([0,t_0],E)}:=\sup_{0\leq t\leq t_0}\E^{-\varrho t}\left( \left\|v(t)\right\|_{E}+ \int_0^t\frac{\left\|v(t)-v(\tau)\right\|_{E}}{(t-\tau)^{\eta+1}}\D\tau\right) <\infty\; ,\]
where $\varrho\geq 1$ is a parameter, \cite{HZ09b}. This standard technique had been used before in \cite{MN03} and \cite{NR02}.

\begin{proposition}\label{P:main}
Assume $0<\alpha,\beta,\gamma,\delta<1$, $\alpha<\gamma<1-\alpha$, $\delta\geq \beta$ and
\[2\gamma+\left(\delta\vee \frac{d_S}{2}\right)<2-2\alpha-\left(\beta\vee \frac{d_S}{2}\right) \; .\]
Let $z\in C^{1-\alpha}([0,t_0], H^{-\beta}(\mu))$ and let $G\in C^2(\mathbb{R})$ with $G(0)=0$ and bounded and Lipschitz second derivative $G''$. Suppose that $R>0$ is given.
Then
\begin{equation}\label{E:invarianceestimate}
\left\|\int_0^{\cdot}T(\cdot-s)G(u(s))\D z(s)\right\|^{(\varrho)}_{W^\gamma([0,t_0], H_{\infty}^\delta(\mu))} \leq C(\varrho)\left(1 +\left\|u\right\|^{(\varrho)}_{W^\gamma([0,t_0],H_{\infty}^\delta(\mu))}\right)\; ,
\end{equation}
$u\in W^\gamma([0,t_0],H_{\infty}^\delta(\mu))$, where $C(\varrho)>0$ tends to zero as $\varrho$ goes to infinity. For sufficiently large $\varrho_0\geq 1$
the closed ball
\[B^{(\varrho_0)}(0,R)=\left\lbrace v\in W^\gamma([0,t_0],H_{\infty}^\delta(\mu)): \left\|v\right\|^{(\varrho_0)}_{W^\eta([0,t_0],H_{\infty}^\sigma(\mu))}\leq R\right\rbrace \]
is mapped into itself and for $\varrho\geq \varrho_0$ large enough,
\begin{multline}
\left\|\int_0^{\cdot}T(\cdot-s)G(u(s))\D z(s)-\int_0^\cdot T(\cdot-s)G(v(s))\D z(s)
\right\|_{W^\gamma([0,t_0],H_{\infty}^\delta(\mu))}^{(\varrho)}\notag\\
\leq C(\varrho)\left\|u-v\right\|^{(\varrho)}_{W^\gamma([0,t_0],H_{\infty}^\delta(\mu))}\;  ,\notag
\end{multline}
$u,v\in B^{(\varrho_0)}(0,R)$.
\end{proposition}

\section{Parabolic problems on metric measure spaces}\label{S:nonlinear}

One of the classes of problems we are interested in are Cauchy problems associated with perturbed semilinear equations. Formulated in a general and abstract way they read
\begin{equation}\label{E:firstex}
\begin{cases}
\frac{\partial u}{\partial t}(t,x)=-\mathrm{A} u(t,x)+F(u)(t,x)+G(u)\cdot\dot{z}(t,x)\ \ ,\ t\in (0,t_0),\ \ x\in X\\
u(0,x)=u_0(x)\; ,
\end{cases}
\end{equation}
where $-\mathrm{A}$ is the $L_2(\mu)$-generator of a strongly continuous symmetric Markovian semigroup $(\mathrm{T}(t))_{t\geq 0}$ on $L_2(\mu)$ as in Section \ref{S:Prelim}, $F$ and $G$ are (generally nonlinear) functions on $\mathbb{R}$ and $\dot{z}$ denotes a space-time perturbation that may be seen as a formal space-time derivative of a non-differentiable deterministic function $z$ on $(0,t_0]\times X$. It is possible to study (\ref{E:firstex}) on general \emph{$\sigma$-finite measure spaces} $(X,\mathcal{X},\mu)$. We will focus on cases with $u(t,x)$ real-valued. In a similar manner one can consider equations with $\mathbb{R}^k$-valued $u(t,x)$, see Example \ref{Ex:parabolic} (ii) and \cite{HZ09b}.

As mentioned before, we investigate the existence of function solutions to these equations. More precisely, we aim at results that confirm the existence and uniqueness of a vector valued function $t\mapsto u(t)$ that solves problem (\ref{E:firstex}) in an evolution sense and takes its values in a space of (equivalence classes of) locally integrable functions on $X$. This is to be distinguished from distribution solutions which are also commonly used to study stochastic partial differential equations.

A function $u$ on $(0,t_0]\times X$ is called a \emph{mild solution} to  (\ref{E:firstex}) if seen as vector-valued function $u(t):=u(t,\cdot)$, it satisfies
\begin{multline}\label{E:mildsol}
u(t)=\mathrm{T}(t)u_0+\int_0^t\mathrm{T}(t-s)F(u(s))\D s +\int_0^t \mathrm{T}(t-s)G(u(s))\D z(s) ,\ \ t\in (0,t_0)\; .
\end{multline}
If for any fixed $t\in [0,t_0]$, $u(t)$ determines a locally integrable function on $(X,\mathcal{X},\mu)$, we call $u$ a \emph{function solution}.
The last term in (\ref{E:mildsol}) is the integral operator as defined in (\ref{E:integral}). It realizes a temporal differentiation of $z$ by means of fractional calculus, a spatial differentiation is hidden in the fact that for fixed time $s$, $z(s)$ is an element of the dual of an appropriate potential space.

The following result is true without any further hypotheses. A proof is given in \cite{HZ12}, its main ingredient is Proposition \ref{P:main}, which allows to use a contraction principle.

\begin{theorem}\label{T:main} Assume $(X,\mathcal{X},\mu)$ is a $\sigma$-finite measure space and $t_0>0$. Let $-\mathrm{A}$ be the generator of a strongly continuous symmetric Markovian semigroup $(\mathrm{T}(t))_{t\geq 0}$ on $L_2(\mu)$ which is ultracontractive with spectral dimension  $d_S>0$.\\
Suppose $0<\alpha,\beta,\gamma,\delta,\varepsilon<1$ and $z\in C^{1-\alpha}([0,t_0], H^{-\beta}(\mu))$. Let $F\in C^1(\mathbb{R})$, $F(0)=0$, have a  bounded Lipschitz derivative $F'$ and $G\in C^2(\mathbb{R})$, $G(0)=0$, have a bounded Lipschitz second derivative $G''$. Assume $f\in H^{2\gamma+\delta+\varepsilon}(\mu)$. If $\alpha<\gamma<1-\alpha$, $\delta\geq \beta$ and
\begin{equation}\label{E:maincond}
2\gamma+\left(\delta\vee \frac{d_S}{2}\right)<2-2\alpha-\left(\beta\vee \frac{d_S}{2}\right) \; .
\end{equation}
Then problem (\ref{E:firstex}) has a unique mild solution (\ref{E:mildsol}) in $W^\gamma([0,t_0], H_{\infty}^{\delta}(\mu))$, which means in particular that the solution is a function.
\end{theorem}

In many cases more structural knowledge about the space $X$ and the semigroup $(\mathrm{T}(t))_{t\geq 0}$ is available. For instance, $X$ may be a \emph{metric measure space} and $(\mathrm{T}(t))_{t\geq 0}$ may possess transition densities that satisfy some typical estimates, \cite{HuZ05, HuZ09}. Under the following assumptions we can improve our results.
\begin{assumption}\label{A:MMS}
$(X,d)$ is a locally compact separable metric space, $\mathcal{X}=\mathcal{B}(X)$ the Borel-$\sigma$-field on $X$ and $\mu$ a Radon measure on $(X,d)$.
\end{assumption}

\begin{assumption}\label{A:HKE}
The semigroup $(\mathrm{T}(t))_{t\geq 0}$ admits transition densities $p(t,x,y)$, that is
\[\mathrm{T}(t)u(x)=\int_X p(t,x,y)u(y)\mu(\D y)\; ,\]
and the $p(t,x,y)$ satisfy bounds of the form
\[t^{-d_f/w}\Phi_1\left(t^{-1/w}d(x,y)\right)\leq p(t,x,y)\leq t^{-d_f/w}\Phi_2\left(t^{-1/w}d(x,y)\right)\]
for any $(x,y)\in X\times X$ and $t\in (0,R_0)$, with bounded decreasing functions $\Phi_i$ on $[0,\infty)$. Here $R_0>0$ is a fixed number, $d_f$ is the \emph{Hausdorff-Dimension} of $(X,d)$ and $w\geq 2$ satisfies $d_S={2d_f}/{w}$. For a given number $\beta>0$ we further assume the validity of the integral condition
\[\int_0^\infty s^{d_f+\beta/2-1}\Phi_2(s)\D s<\infty\; .\]
\end{assumption}

Under these circumstances we get the following improved result.

\begin{theorem}\label{T:main2} Let $F$ and $G$ be as in Theorem \ref{T:main}. Suppose $0<\alpha,\beta,\gamma,\delta,\varepsilon<1$ and Assumptions \ref{A:MMS} and \ref{A:HKE} are satisfied. Assume  $\alpha<\gamma<1-\alpha$ and $0<\beta<\delta<{d_S}/{2}$. If $f\in H^{2\gamma+\delta+\varepsilon}(\mu)$ and $z\in C^{1-\alpha}([0,t_0], H_{q}^{-\beta}(\mu))$ for $q={d_S}/{\delta}$ and
\begin{equation}\label{E:maincond2}
2\gamma+\frac{d_S}{2}<2-2\alpha-\beta\; ,
\end{equation}
then problem (\ref{E:firstex}) has a unique mild solution (\ref{E:mildsol}) in $W^\gamma([0,t_0], H_{\infty}^{\delta}(\mu))$, which means in particular that the solution is a function.
\end{theorem}

Theorem \ref{T:main2} is proved in a similar way as Theorem \ref{T:main} by verifying the contractivity of the integral operator and applying a contraction argument. In particular, analogs of Lemmas \ref{L:WeylMarchaud} and \ref{L:existence} can be used. The only news is the following improved product estimate that replaces the former (\ref{E:basicproduct}).
\begin{proposition}
Let $0<\beta<\delta< {d_S}/{2}\wedge 1$ and $p={d_S}/{(d_S-\delta)}$. Let the semigroup be ultracontractive with spectral dimension $d_S>0$ and let Assumptions \ref{A:MMS} and \ref{A:HKE} be satisfied. Then we have
\[\left\|uv\right\|_{H_p^\beta(\mu)}\leq c\:\left\|u\right\|_\delta\left\|v\right\|_\delta\]
for any $u,v\in H^\delta(\mu)$.
\end{proposition}

Theorem \ref{T:main2} requires $d_S<4$. For symmetric diffusion semigroups on $\mathbb{R}^n$ we have $d_S=n$, hence need $n\leq 3$. This is typical, because to deal with the nonlinear transformations $F$ and $G$ we need the solution to be $L_\infty(\mu)$-bounded, but only in low dimensions the singularity of the semigroup at zero is small enough to provide $L_\infty(\mu)$-bounds. The special case of linear $F$ and $G$ allows to remove this restrictive condition.

\begin{theorem}\label{T:main3} Let $F$ and $G$ be linear. Suppose $0<\alpha,\beta,\gamma,\delta,\varepsilon<1$ and Assumptions \ref{A:MMS} and \ref{A:HKE} are satisfied. Assume  $\alpha<\gamma<1-\alpha$ and
$0<\beta<\delta<{d_S}/{2}$. If $f\in H^{2\gamma+\delta+\varepsilon}(\mu)$ and $z\in C^{1-\alpha}([0,t_0], H_{q}^{-\beta}(\mu))$ for $q={d_S}/{\delta}$ and
\begin{equation}\label{E:maincond3}
2\gamma+\delta<2-2\alpha-\beta\; ,
\end{equation}
then problem (\ref{E:firstex}) has a unique mild solution (\ref{E:mildsol}) in $W^\gamma([0,t_0], H^{\delta}(\mu))$, which means in particular that the solution is a function.
\end{theorem}

\begin{examples}\label{Ex:parabolic}
To consider some examples of stochastic partial differential equations based on (\ref{E:firstex}), let $0<H,K<1$ and consider the spatially isotropic fractional Brownian sheet $B^{H,K}$ on $[0,t_0]\times\mathbb{R}^n$ with Hurst indices $H$ and $K$ (see \cite{ALP02}), that is, the centered real valued Gaussian random field $B^{H,K}$ on $[0,t_0]\times\mathbb{R}^n$ over a probability space $(\Omega,\mathcal{F},\mathbb{P})$ such that for any $0\leq s<t\leq t_0$ and $x,y\in\mathbb{R}^n$, 
\[\mathbb{E}\left[B^{H,K}(t,y)-B^{H,K}(t,x)-B^{H,K}(s,y)+B^{H,K}(s,x)\right]^2 =c_{H,K}(t-s)^{2H}|x-y|^{2K}\; ,\]
where $|x|$ denotes the Euclidean norm of $x\in\mathbb{R}^n$. It is not difficult to see that for $\mathbb{P}$-a.e. $\omega\in\Omega$ and any
$0<\gamma<H$, $0<\sigma<K$ and $1<q<\infty$, the realization $B^{H,K}(\omega)$ is an element of $C^\gamma([0,t_0], H_q^\sigma(\mathbb{R}^n))$. The components of its distributional (spatial) gradient $\nabla B^{H,K}(\omega)$ are elements of $C^\gamma([0,t_0], H_q^{\sigma-1}(\mathbb{R}^n))$.
\begin{enumerate}
\item[(i)] Let $n=1$, $X=(0,1)$, let $(\mathrm{T}(t))_{t\geq 0}$ be the Dirichlet heat semigroup on $(0,1)$ and $\varDelta$ the Dirichlet Laplacian. Consider the one-dimensional semilinear heat equation on $(0,t_0)\times (0,1)$ driven by a fractional Brownian sheet $B^{H,K}$,
\[\frac{\partial u}{\partial t}(t,x)=\varDelta u(t,x)+ F(u(t,x))+G(u(t,x))\cdot\frac{\partial^2 B^{H,K}}{\partial t\partial x}\; .\]
It has a unique function solution if $1/2<H<1$ and $2H+K>2$.
\item[(ii)] In \cite{HZ09b} we have considered boundary initial value problems on smooth bounded domains $D\subset\mathbb{R}^n$ associated with parabolic equations of type
\[\frac{\partial u}{\partial t}(t,x)=-\mathrm{A}u(t,x)+F(u(t,x))+\left\langle G(u),\frac{\partial}{\partial t}\nabla V\right\rangle_{\mathbb{R}^n}(t,x)\; .\]
Here $V$ is a real valued noise potential, $G$ is an $\mathbb{R}^n$-valued nonlinearity on $\mathbb{R}$ and $\left\langle \cdot, \cdot\right\rangle_{\mathbb{R}^n}$ denotes the scalar product in $\mathbb{R}^n$.
\end{enumerate}
\end{examples}

\section{Transport equations on domains}

In this section we consider transport-diffusion equations of form
\begin{equation}\label{E:transport}
\begin{cases}
\frac{\partial u}{\partial t}(t,x)=\varDelta u(t,x)+\left\langle \nabla u,\nabla z\right\rangle_{\mathbb{R}^n}(t,x), & t\in (0,t_0],\  x\in D\\
u(t,x)=0, & t\in (0,t_0],\  x\in\partial D\\
u(0,x)=u_0(x), & x\in D\; ,
\end{cases}
\end{equation}
where $D\subset\mathbb{R}^n$ is a smooth bounded domain and $z$ is a non-differentiable function on $\mathbb{R}^n$. Here $\varDelta$ denotes the Dirichlet Laplacian for $D$ and the gradient $\nabla$ is interpreted in distributional sense. As before, $\left\langle\cdot,\cdot\right\rangle_{\mathbb{R}^n}$ denotes the scalar product in $\mathbb{R}^n$. In this model $z$ is viewed as a temporally constant perturbation.

Problems of type (\ref{E:transport}) have been considered by the second named author in \cite{I}. Again we are interested in the existence, uniqueness and regularity of function solutions. Now $(\mathrm{T}(t))_{t\geq 0}$ will denote the Dirichlet heat semigroup $(\mathrm{T}^D(t))_{t\geq 0}$ for the domain $D$.

A function $u$ on $(0,t_0]\times D$ is called a \emph{mild solution} to (\ref{E:transport}) if seen as a vector valued function $u(t):=u(t,\cdot)$ it satisfies
\[u(t)=\mathrm{T}(t)u_0+\int_0^t \mathrm{T}(t-r)\left\langle \nabla u(r),\nabla z\right\rangle_{\mathbb{R}^n} \D r, \ \ t\in (0,t_0]\; ,\]
and a \emph{function solution} if in addition $u(t)$ is a locally integrable function on $D$ for any $t$. By the following multiplication property the right-hand side in this definition admits a reasonable interpretation.

\begin{lemma}
Let $w\in \widetilde{H}^{1+\delta}_p(D)$, $z\in H_q^{1-\beta}(\mathbb{R}^n)$ with $1<p,q<\infty$, $q>p\vee {n}/{\delta}$, $0<\beta<1/2$ and $\beta<\delta$. Then $\left\langle \nabla w, \nabla z\right\rangle_{\mathbb{R}^n}$ is a member of $\widetilde{H}_p^{-\beta}(D)$ and
\[\left\|\left\langle \nabla w,\nabla z\right\rangle_{\mathbb{R}^n}\right\|_{H_p^{-\beta}(\mathbb{R}^n)}\leq c\left\|w\right\|_{H^{1+\delta}_p(\mathbb{R}^n)} \left\|z\right\|_{H^{1-\beta}_q(\mathbb{R}^n)}\; .\]
\end{lemma}

The main result of \cite{I} reads as follows.

\begin{theorem}\label{T:maintransport}
Let $t_0>0$ and $0<\beta<\delta<1/2$ and $0<2\gamma<1-\beta-\delta$. Let $z\in H_q^{1-\beta}$ for some $q>2\vee {d}/{\delta}$.
Then for any $u_0\in \widetilde{H}^{1+\delta+2\gamma}(D)$ there exists a unique mild solution $u\in C^\gamma([0,t_0], \widetilde{H}^{1+\delta}(D))$ to (\ref{E:transport}), which means in particular that the solution is a function.
\end{theorem}

The theorem follows by fixed point arguments and the following contractivity result. Similarly as before it is formulated in terms of equivalent norms. For $\varrho\geq 1$ we equip the space $C^\eta([0,t_0],E)$ of $\eta$-H\"older continuous $E$-valued functions $v$ on $[0, t_0]$ with the equivalent norm
\[\left\|v\right\|_{C^\eta([0,t_0],E)}^{(\varrho)}:=\sup_{0\leq t\leq t_0} \E^{-\varrho t}\left(\left\|v(t)\right\|_{E}+\sup_{0\leq\tau<t} \frac{\left\|v(t)-v(\tau)\right\|_{E}}{(t-\tau)^\gamma}\right)\; .\]

\begin{proposition}
$0<\beta<\delta<1/2$ and $z\in H_q^{1-\beta}(\mathbb{R}^d)$ for some $q>2\vee {d}/{\delta}$. Then for any $\gamma$ with $0<2\gamma<1-\beta-\delta$ we have
\[\left\|\int_0^\cdot \mathrm{T}(\cdot-s)\left\langle \nabla u(s),\nabla z\right\rangle_{\mathbb{R}^n} ds\right\|_{C^\gamma([0,t_0],\widetilde{H}^{1+\delta}(D))}^{(\varrho)}
\leq C(\varrho)\left\|u\right\|_{C^\gamma([0,t_0],\widetilde{H}^{1+\delta}(D))}^{(\varrho)}\]
for any $u\in C^\gamma([0,t_0],\widetilde{H}^{1+\delta}(D))$ where $C(\varrho)$ tends to zero as $\varrho$ goes to infinity.
\end{proposition}

\begin{examples}
\begin{enumerate}
\item[(i)] If for instance $B^H$ is a fractional Brownian field on $\mathbb{R}^n$ with Hurst parameter $1/2<H<1$, that is a real valued centered Gaussian random field on $\mathbb{R}^n$ with
\[\mathbb{E}\left[B^H(x)-B^H(y)\right]^2=c_H|x-y|^{2H}\; ,\]
then we may consider a typical realization $B^H(\omega)$ in place of $z$ to obtain results for stochastic transport equations
with fractal noise,
\[\frac{\partial u}{\partial t}(t,x)=\varDelta u(t,x)+\left\langle \nabla u,\nabla B^H\right\rangle_{\mathbb{R}^n} (t,x)\; .\]
\item[(ii)]
Combining the above with the results of the preceding section we can investigate a more general form of transport-diffusion equation,
\[\begin{cases}
\frac{\partial u}{\partial t}(t,x)=\varDelta u(t,x)+\left\langle \nabla u,\nabla z\right\rangle_{\mathbb{R}^n}(t,x)+\left\langle F,\frac{\partial}{\partial t}\nabla V\right\rangle_{\mathbb{R}^n}, & t\in (0,t_0],\  x\in D\\
u(t,x)=0, & t\in (0,t_0],\  x\in\partial D\\
u(0,x)=u_0(x), & x\in D\; ,
\end{cases}\]
where $z$ is as above, $F$ is a vector in $\mathbb{R}^n$ and $V=V(t,x)$ is a non-differentiable noise that may vary in space and time.
\end{enumerate}
\end{examples}

\section{Some remarks on Burgers system}

We finish our survey with some brief look at a \emph{Burgers type equation}, \cite{B74}. On $(0,t_0)\times\mathbb{R}^n$, $n\geq 1$, consider the equation
\begin{equation}\label{E:BurgersPDE}
\frac{\partial u}{\partial t} = \varDelta u -\left\langle u,\nabla\right\rangle_{\mathbb{R}^n} u + \frac{\partial}{\partial t}\nabla B\ \
\end{equation}
with some deterministic initial condition $u(0)=u_0$. Here $\left\langle \cdot,\cdot \right\rangle_{\mathbb{R}^n}$ denotes the Euclidean scalar product in $\mathbb{R}^n$ and $B=B(t,x)$ is a \emph{fractional Brownian sheet} on $[0,t_0]\times\mathbb{R}^n$ over some probability space $(\Omega,\mathcal{F},\mathbb{P})$. Equation (\ref{E:BurgersPDE}) is already a stochastic differential equation and the solution method considered by the first named author in \cite{H09} is not pathwise. However, it is made up from techniques very similar to those used in the preceding sections. Note that a solution $u$ will be vector valued, i.e. $u(t,x)\in\mathbb{R}^n$ for fixed $t$ and $x$.

We make (\ref{E:BurgersPDE}) rigorous by defining weak and mild solutions. A process $u=u(t)$ is said to be a \emph{distributional solution} to equation (\ref{E:BurgersPDE}) on $(0,t_0)$ with initial condition $u_0$ if for any test function $\varphi\in \mathcal{D}(\mathbb{R}^n,\mathbb{R}^n)$ and any $t\in (0,t_0)$, we $\mathbb{P}$-a.s. have
\begin{equation}\label{E:solBurgers}
(u(t),\varphi)=(u_0,\varphi)+\int_0^t(\varDelta u(s),\varphi)\D s-\int_0^t(\left\langle u(s),\nabla\right\rangle_{\mathbb{R}^n} u(s),\varphi)\D s+ (\nabla B(t),\varphi)\; .
\end{equation}
Here $\varDelta$ denotes the matrix Laplacian, that is the $n\times n$ diagonal matrix with the usual Laplacian on the diagonal. $B(t)=B(t,\cdot)$ can be seen as a process taking values in a Sobolev space and the components $\partial B(t)/\partial x_i$ of $\nabla B(t)=(\partial B(t)/\partial x_1, ... ,\partial B(t)/\partial x_n)$ are defined in the sense of distributions. As the name indicates, $u$ is considered as a distribution-valued function. However, it turns out that it also provides a function solution. A process $u=u(t)$ is called a \emph{mild solution} to (\ref{E:BurgersPDE}) on $(0,t_0)$ with initial condition $u_0$ if for any $t\in (0,t_0)$, $\mathbb{P}$-a.s.
\begin{equation}\label{E:solBurgersmild}
u(t)=\mathrm{T}(t)u_0-\int_0^t \mathrm{T}(t-s) \left\langle u(s),\nabla\right\rangle_{\mathbb{R}^n} u(s)\D s+\int_0^t \mathrm{T}(t-s)\D(\nabla B)(s)\; .
\end{equation}
Here $(\mathrm{T}(t))_{\geq 0}$ denotes the Brownian semigroup on $\mathbb{R}^n$, respectively its matrix version. The stochastic integral in (\ref{E:solBurgersmild}) is of Wiener type, defined as the $L_p(\Omega)$-limit
\begin{multline}
\int_0^t \mathrm{T}(t-s)\D (\nabla B)(s):=\lim_{\varepsilon\to 0}\varepsilon\int_0^1r^{\varepsilon-1}\int_0^t\mathrm{T}(t-s)\frac{\nabla B(s+r)-\nabla B(s)}{r} \D s\D r\; ,
\end{multline}
$p>2$, of random variables taking values in a certain weighted Sobolev space. It is closely related to the pathwise integrals (\ref{E:vectorforward}) and (\ref{E:integral}). We call $u$ as in (\ref{E:solBurgersmild}) a \emph{function solution} to (\ref{E:BurgersPDE}) if any $u(t)$, $t\in (0,t_0]$, is a locally integrable function on $\mathbb{R}^n$.

We follow the standard approach to Burgers equation and employ a stochastic variant of the \emph{Cole-Hopf transformation}. First we consider a related \emph{stochastic heat equation},
\begin{equation}\label{E:halfnoiseeq}
\frac{\partial w}{\partial t}(t,x)= \varDelta w(t,x) + w(t,x)\cdot\frac{\partial}{\partial t} B(t,x) \ \ ,\ t\in (0,T)\ ,\ x\in\mathbb{R}^n\; ,
\end{equation}
with some initial condition $w_0$. The term $\partial B /\partial t$ is a \emph{half-noise}, similar as in \cite{BCJL94}: For fixed $x\in \mathbb{R}^n$ and up to a constant, $t\mapsto B(t,x)$ behaves like a one-dimensional fractional Brownian motion with Hurst parameter $H$ and in (\ref{E:halfnoiseeq}) we consider its formal time derivative $t\mapsto \partial B(t,x) /\partial t $. We say the random process $w:(0,t_0)\times\mathbb{R}^n\times\Omega\to \mathbb{R}$ is a \emph{pointwise mild solution} to (\ref{E:halfnoiseeq}) if for fixed  $t\in (0,t_0)$ and $x\in\mathbb{R}^n$, we have
\begin{equation}\label{E:halfnoisesol}
w(t,x)=\mathrm{T}(t)w_0(x)+\lim_{\varepsilon\to 0}\varepsilon\int_0^1r^{\varepsilon-1}\int_0^t \mathrm{T}(t-s)w(s,x)\frac{B(s+r,x)-B(s,x)}{r}\D s\D r\; .
\end{equation}
The limit and the equality (\ref{E:halfnoisesol}) are considered in $L_p(\Omega)$, $p>1$. Our results are as follows:

\begin{theorem}\label{T:dist}
Let $t_0>0$, $0<K\leq 1/2$ and $2<2H+K$. Suppose that $u_0$ is of form $u_0(x)=-\nabla U_0(x)$, where $U_0$ is a real-valued function on $\mathbb{R}^n$ such that
\[|U_0(x)|\leq b(1+|x|^{\gamma})\;,\ \ \left|\frac{\partial U_0}{\partial x_i}(x)\right|\leq\exp(b(1+|x|^{\gamma}))\; ,\ \  i=1,...,n\]
for some $b>0$, $2K\leq\gamma\leq 1$, and any $x\in\mathbb{R}^n$. Then there is a $\mathbb{P}$-a.s. strictly positive pointwise mild solution $w$ to (\ref{E:halfnoiseeq}) with $w_0=\exp(U_0/2)$ and $u:=\nabla\log w$ is a distributional solution (\ref{E:solBurgers}) to (\ref{E:BurgersPDE}).
The process $u$ is also a function solution to (\ref{E:halfnoiseeq}).
\end{theorem}

Note that our hypothesis implies $H>3/4$. For further details we refer to \cite{H09}.

%\section*{Appendix}
%\addcontentsline{toc}{section}{Appendix}

%%%%%%%%%%%%%%%%%%%%%%%% references.tex %%%%%%%%%%%%%%%%%%%%%%%%%%%%%%

\end{document}